\documentclass{amsart}
\usepackage{amssymb}
\usepackage{amsfonts}
\usepackage{amsmath}
\usepackage[legalpaper,bookmarks=true,colorlinks=true,linkcolor=blue,citecolor=blue]{hyperref}
\usepackage{graphicx}%
\usepackage{fancyhdr}
\usepackage{color}
\usepackage[mathlines]{lineno}
\usepackage{lscape}
\usepackage{epsfig}
\usepackage{natbib}
\usepackage{geometry}
\usepackage{tgbonum}
\usepackage[]{listings}

\fontfamily{qcr}\selectfont

\newtheorem{theorem}{Theorem}
\theoremstyle{plain}

\newtheorem{lemma}{Lemma}

\numberwithin{equation}{section}

\newcommand{\Bin}{\bigskip \noindent}

\newcommand{\Ni}{\noindent}

\begin{document}
\Large
\title[General asymptotic representations of indexes using the fep-lorep]{General asymptotic representations of indexes based on the functional empirical process and the residual functional empirical process and applications}
\maketitle
\noindent \author{Gane Samb Lo $^{\dag}$, Tchilabalo Abozou Kpanzou $^{\dag\dag}$, Gandasor Bonyiri Onesiphore $^{\dag\dag\dag}$}

\bigskip
\noindent $^{\dag}$ Gane Samb Lo.\\
LERSTAD, Gaston Berger University, Saint-Louis, S\'en\'egal (main affiliation).\newline
LSTA, Pierre and Marie Curie University, Paris VI, France.\newline
AUST - African University of Sciences and Technology, Abuja, Nigeria\\
gane-samb.lo@edu.ugb.sn, gslo@aust.edu.ng, ganesamblo@ganesamblo.net\\
Permanent address : 1178 Evanston Dr NW T3P 0J9, Calgary, Alberta, Canada.\\

\noindent $^{\dag\dag}$ Dr Tchilabalo Abozou Kpanzou
University of Kara, Kara, Togo\\
Affiliated to LERSTAD, Gaston Berger University, Saint-Louis, Senegal\\
Emails: t.kpanzou@univkara.net; kpanzout@gmail.com\\

\noindent $^{\dag\dag\dag}$ Da Gandasor Bonyiri Onesiphore\\
LERSTAD, Gaston Berger University (UGB), Saint-Louis, Senegal.\\
IMHOTEP International Mathematical Centre (IMHO-IMC: imhotepsciences.org)\\

\newpage \label{fsgp2gslo}
\Ni \textbf{Full Abstract}. The objective of this paper is to establish a general asymptotic representation (\textit{GAR}) for a wide range of statistics, employing two fundamental processes: the functional empirical process (\textit{fep}) and the residual functional empirical process introduced by Lo and Sall (2010a, 2010b), denoted as \textit{lrfep}. The functional empirical process (\textit{fep}) is defined as follows:

$$
\mathbb{G}_n(h)=\frac{1}{\sqrt{n}} \sum_{j=1}^{n} \{h(X_j)-\mathbb{E}h(X_j)\},
$$

\Bin [where $X$, $X_1$, $\cdots$, $X_n$ is a sample from a random $d$-vectors $X$ of size $(n+1)$ with  $n\geq  1$ and $h$ is a measurable function defined on $\mathbb{R}^d$ such that $\mathbb{E}h(X)^2<+\infty$]. It is a powerful tool for deriving asymptotic laws. An earlier and simpler version of this paper focused on the application of the (\textit{fep}) to statistics $J_n$ that can be turned into an asymptotic algebraic expression of empirical functions of the form

$$
J_n=\mathbb{E}h(X) + n^{-1/2} \mathbb{G}_n(h) + o_{\mathbb{P}}(n^{-1/2}).  \ \ \ \textit{SGAR}
$$

\Bin  However, not all statistics, in particular welfare indexes, conform to this form. In many scenarios, functions of the order statistics $X_{1,n}\leq$, $\cdots$, $\leq X_{n,n}$ are involved, resulting in $L$-statistics. In such cases, the (\textit{fep}) can still be utilized, but in combination with the related residual functional empirical process introduced by Lo and Sall (2010a, 2010b). This combination leads to general asymptotic representations (GAR) for a wide range of statistical indexes

$$
J_n=\mathbb{E}h(X) + n^{-1/2} \biggr(\mathbb{G}_n(h) + \int_{0}^{1} \mathbb{G}_n(\tilde{f}_s) \ell(s) \ ds + o_{\mathbb{P}}(1)\biggr), \ \ \textit{FGAR}
$$

\Bin where $\tilde{f}_s=1_{]-\infty, F^{-1}(s)]}$ is the indicator function of the interval $]-\infty, F^{-1}(s)]$ of $\mathbb{R}$, $h \in L^2(\mathbb{P}_X)$ and $\ell(\circ)$ is a measurable function of $s \in (0,1)$. Such representations, when associated with copulas, provide a robust framework for comparing indices  over the time or across different areas. The comprehensive theory is presented alongside explicit examples, facilitating the utilization by  researchers.\\
\newpage

\noindent \textbf{R\'{e}sum\'{e}.} (French abstract). L'objectif de cet article est d'établir une représentation asymptotique générale (\textit{RAG}) pour un large éventail de statistiques, en utilisant deux processus : le processus empirique fonctionnel (\textit{fep}) et le processus empirique résiduel introduit par Lo et Sall (2010a, 2010b), abrégé sous le nom de \textit{lrfep}. Le processus empirique fonctionnel (\textit{fep}) est défini comme suit :

$$
\mathbb{G}_n(h)=\frac{1}{\sqrt{n}} \sum_{j=1}^{n} \{h(X_j)-\mathbb{E}h(X_j)\},
$$

\Bin [où $X$, $X_1$, $\cdots$, $X_n$ représentent un échantillon d'un vecteur aléatoire (de longueur $d$) de taille 
$(n+1)$ with $n\geq 1$, et $h$ est une fonction mesurable définie sur $\mathbb{R}^d$ telle que 
$\mathbb{E}h(X)^2<+\infty$]. Ce processus est un outil puissant pour déduire les propriétés asymptotiques. Une version antérieure et plus simple de cet article se concentrait sur l'application du \textit{fep} aux statistiques $J_n$ qui peuvent être transformées en une expression algébrique asymptotique des fonctions empiriques sous la forme :

$$
J_n=\mathbb{E}h(X) + n^{-1/2} \mathbb{G}_n(h) + o_{\mathbb{P}}(n^{-1/2}). \ \ \textit{SGAR}
$$

\Bin Cependant, toutes les statistiques, en particulier les indices de bien-être, ne suivent pas cette forme. Dans de nombreuses situations, des fonctions des statistiques d'ordre $X_{1,n}\leq$, $\cdots$, $\leq X_{n,n}$ sont impliquées, ce qui donne naissance aux statistiques $L$. Dans de tels cas, le \textit{fep} peut toujours être utilisé, mais en combinaison avec le processus empirique résiduel associé introduit par Lo et Sall (2010a, 2010b). Cette combinaison conduit à une représentation asymptotique générale (GAR) pour un large éventail d'indices statistiques :

$$
J_n=\mathbb{E}h(X) + n^{-1/2} \biggr(\mathbb{G}_n(h) + \int_{0}^{1} \mathbb{G}_n(\tilde{f}_s) \ell(s) \ ds + o_{\mathbb{P}}(1)\biggr), \ \ \textit{FGAR}
$$

\Bin o\`u $\tilde{f}_s=1_{]-\infty, F^{-1}(s)]}$ est la fonction indicatrice de l'intervalle  $]-\infty, F^{-1}(s)]$ de $\mathbb{R}$, $h \in L^2(\mathbb{P}_X)$, et $\ell(\circ)$ est une fonction mesurable de $s \in (0,1)$. De telles représentations, lorsqu'elles sont combinées avec les copules, fournissent un cadre solide pour la comparaison des indices au fil du temps ou entre différentes zones. La théorie complète est présentée avec des exemples explicites, facilitant son utilisation par les chercheurs.\\

\newpage
\Bin \textbf{Presentations of authors}.\\

\Ni \textbf{Gane Samb Lo}, Ph.D., is a retired full professor from Université Gaston Berger (2023), Saint-Louis, SENEGAL. He is the founder and The Probability and Statistics Chair holder at: Imhotep International Mathematical Center (imho-imc), https://imhotepsciences.org.\\

\Ni \textbf{Tchilabalo Abozou Kpanzou}, Ph.D.,  is a professor of Statistics and Mathematics at Universite of Kara, Togo. He is a lead researcher at: Imhotep International Mathematical Center (imh-imc), https://imhotepsciences.org.\\

\Ni \textbf{Gandasor Bonyiri Onesiphore Da}, M.Sc, is preparing a Ph.D. degree at: Doctoral Studies, LERSTAD, Université Gaston Berger, Saint-Louis, SENEGAL, and in the Imhotep International Mathematical Center (imh-imc), https://imhotepsciences.org, under the supervision of the previous two co-authors.\\

\newpage
\section{Introduction} \label{sec1}

\noindent  There is a huge number of statistics and indexes, that are used in Mathematical Statistics and in many other disciplines, whose asymptotic behaviour can be studied in a unified approach based on the functional empirical  processes. For instance, indexes of poverty have been extensively used in Ecometrics. A detailed survey of them is available in \cite{zeng}, including simple ones as the Foster, Greer and Thorbecke index (See \cite{fgt}) and  most elaborated ones such as the statistics of \cite{kakwani} and 
\cite{sen} (Economics Nobel prize winner). The first attempt of having a common frame for the asymptotic theory of those statistics was done in \cite{logpi}. However, the approach in that paper uses the real-valued empirical processes. Later, it was realized that the rightest approach would rely on function-value empirical process, i.e., the functional empirical process (\textit{fep}). The combination of the \textit{fep} with the so-called residual empirical process - that has been later turned into a functional form named \textit{residual functional empirical process} - started a new and global approach that provides general asymptotic representations (\textit{GAR}) in the \textit{fep} for a huge class of indices arising in many fields. (\textit{GAR})'s have been already provided  in a number of papers, in particular in \cite{mergane2018Gini} for Gini's inequality index, \cite{mergane2018Theil} for the class of Theil indexes. These \textit{GAR}s have proved to be very effective in addressing the decomposability of statistics topic by tackling the lack of decomposability as in \cite{haidaraDecomp}. The main reason of the success of the \textit{fep} is its linearity property that the real-valued empirical process fails to have.\\

\Ni We feel that this method of \textit{GAR}'s in the \textit{fep} is not sufficiently exploited and it is even unknown by the great public of applied researchers. In that sense, the main motivation of this paper is to provide a new and global exposure of this method to researchers in order to show them how to apply it in specific cases and to give explicit examples of usage.\\

\Ni At the same time, a number of applications in a series of papers will follow in the special number of this journal. Another paper on statistical computations and related packages will be published. That sample of cases, more or less numerous, will be the beginning of a general handbook that will be handed to researchers.\\

\Ni What is very amazing is that, despite its power, this technique does not need the heavy weapons of functional weak convergence such as Donsker classes, Vapnik-Chervonenkis theory (see \cite{vcclass}) or entropy numbers, etc., nor the contents of books like \cite{vaart}, \cite{gaenssler}\ or in \cite{pollard} or alike texts. Actually, we only use our tools at the finite-distribution level, where the functional empirical processes can be reasonably easy to use without heavy theories.\\

\Ni Accordingly, we will organize the rest of the paper as follows. In next section (Section \ref{sec_02}), we describe the \textsl{fep}, justify it and consider a case study on a non-trivial example. Namely, we will entirely show how the method works for the plug-in estimator of the linear correlation index between two random variables. In Section \ref{sec_03}, we introduce the \textit{lrfep}, give the conditions ensuring its weak convergence. In Section \ref{sec_04}, we use the Gini's index as a study case for applying the recommendations in the previous section. In Section \ref{sec_05}, we make a final summary of the whole method and give final remarks and recommendations. 

\newpage
\section{The Functional Empirical Process and applications in finding a simple \textit{(GAR)}} \label{sec_02}

\subsection{The functional empirical process \textit{(fep)}} \label{sec_02_subsec_01}

\Ni Let us begin by introducing the \textit{fep} and describe its asymptotic representation.\\

\Ni \textbf{What is the \textit{fep}?}.\\
 
\Ni Let $Z_{1}$, $Z_{2}$, ... be a sequence of independent copies of a random variable $Z$ defined on the same probability space with
values on some metric space $(S,d)$. We define for each $n\geq 1,$ the functional empirical process by 

\begin{equation*}
\mathbb{G}_{n}(f)=\frac{1}{\sqrt{n}}\sum_{i=1}^{n}(f(Z_{i})-\mathbb{E}%
f(Z_{i})),
\end{equation*}

\bigskip \noindent where $f$ is a real and measurable function defined on $\mathbb{R}$ such that

\begin{equation}
\mathbb{V}_{Z}(f)=\int \left( f(x)-\mathbb{P}_{Z}(f)\right)
^{2}dP_{Z}(x)<\infty ,  \label{var}
\end{equation}%
which entails

\begin{equation}
\mathbb{P}_{Z}(\left\vert f\right\vert )=\int \left\vert f(x)\right\vert
dP_{Z}(x)<\infty \text{.}  \label{esp}
\end{equation}

\bigskip \noindent We denote by $\mathcal{F}(S)$ - $\mathcal{F}$ for short -
the class of real-valued measurable functions that are defined on $S$ such
that (\ref{var}) holds. The space $\mathcal{F}$ , when endowed with the
addition and the external multiplication by real scalars, is a linear space.
Next, it is remarkable that $\mathbb{G}_{n}$ is linear on $\mathcal{F}$, that
is for $f$ and $g$ in $\mathcal{F}$ and for $(a,b)\in \mathbb{R}{^{2}}$, we
have

\begin{equation*}
a\mathbb{G}_{n}(f)+b\mathbb{G}_{n}(g)=\mathbb{G}_{n}(af+bg).
\end{equation*}

\bigskip \noindent We have this result.

\begin{lemma} \label{lemma.tool.1}
\bigskip Given the notation above, then for any finite number of elements $%
f_{1},...,f_{k}$ of $\mathcal{F}$, $k\geq 1$, we have

\begin{equation*}
^{t}(\mathbb{G}_{n}(f_{1}),...,\mathbb{G}_{n}(f_{k}))\rightsquigarrow 
\mathcal{N}_{k}(0,\Gamma (f_{i},f_{j})_{1\leq i,j\leq k}),
\end{equation*}

\bigskip \noindent where 
\begin{equation*}
\Gamma (f_{i},f_{j})=\int \left( f_{i}-\mathbb{P}_{Z}(f_{i})\right) \left(
f_{j}-\mathbb{P}_{Z}(f_{j})\right) d\mathbb{P}_{Z}(x),1\leq i,j\leq k.
\end{equation*}
\end{lemma}

\bigskip \noindent \textbf{Proof}. It is enough to use the Cram\'{e}r-Wold Criterion (see for example \cite{billinsgley68}, page 45), that
is to show that for any $a=^{t}(a_{1},...,a_{k})\in \mathbb{R}^{k},$ by
denoting $T_{n}=^{t}(\mathbb{G}_{n}(f_{1}),...,\mathbb{G}_{n}(f_{k})),$ we have $<a,T_{n}>\rightsquigarrow <a,T>$, where $T$ follows the $\mathcal{N}%
_{k}(0,\Gamma (f_{i},f_{j})_{1\leq i,j\leq k})$\ law and $<\circ ,\circ >$
stands for the usual scalar product in $\mathbb{R}^{k}.$ But, by the standard central limit theorem in $\mathbb{R}$, we have%
\begin{equation*}
<a,T_{n}>=\mathbb{G}_{n}\left( \sum\limits_{i=1}^{k}a_{i}f_{i}\right)
\rightsquigarrow N(0,\sigma _{\infty }^{2}),
\end{equation*}

\bigskip \noindent where, for $g=\sum_{1\leq i\leq k}a_{i}f_{i},$%
\begin{equation*}
\sigma _{\infty }^{2}=\int \left( g(x)-\mathbb{P}_{Z}(g)\right) ^{2}dP_{Z}(x)
\end{equation*}

\bigskip \noindent and this easily gives%
\begin{equation*}
\sigma _{\infty }^{2}=\sum\limits_{1\leq i,j\leq k}a_{i}a_{j}\Gamma
(f_{i},f_{j}),
\end{equation*}

\Ni so that $N(0,\sigma _{\infty }^{2})$ is the law of $<a,T>.$ The proof is
finished.

\subsection{How to use the tool?}

\Ni In our daily activities as researchers in asymptotic statistics, we usually need to find the asymptotic laws of statistics on $\mathbb{R}^{k}$. Once we have our sample $Z_{1}$, $Z_{2}$, $...$ as random variables defined on the same probability space with values in $\mathbb{R}
^{k},$ the studied statistic, say $H_{n}$, in a significant number of cases, can be expressed in the following form:

\begin{equation*}
H_{n}=\frac{1}{n}\sum\limits_{i=1}^{k}H(Z_{i})
\end{equation*}

\bigskip \noindent for $H\in \mathcal{F}$. This gives the simplest asymptotic representation,  for $\mu (H)=\mathbb{E}H(Z)$: 

\begin{equation}
H_{n}=\mu (H)+n^{-1/2}\mathbb{G}_{n}(H).  \label{expan}
\end{equation}

\bigskip \noindent From such exact \textit{GAR}'s, the Delta-method is a common tools to proceed to next steps. For that, we need to address the probability asymototic boundedness of $\mathbb{G}_{n}(H)$. Since it weakly converges to a random element, say 
$M(H)$, we get $\left\Vert \mathbb{G}_{n}(H)\right\Vert \rightsquigarrow \left\Vert M(H)\right\Vert $, by the continuous mapping theorem. Next, we get for all $\lambda >0,$ by the assertion of the Portmanteau Theorem concerning open sets,

\begin{equation*}
\lim \sup_{n\rightarrow \infty }P(\left\Vert \mathbb{G}_{n}(H)\right\Vert
>\lambda )\leq P(\left\Vert M(H)\right\Vert >\lambda )
\end{equation*}

\bigskip \noindent and then 
\begin{equation*}
\lim \inf_{\lambda \rightarrow \infty }\lim \sup_{n\rightarrow \infty
}P(\left\Vert \mathbb{G}_{n}(H)\right\Vert >\lambda )\leq \lim \sup
P(\left\Vert M(H)\right\Vert >\lambda )=0.
\end{equation*}

\bigskip \noindent From this, we use the big $O_{\mathbb{P}}$ notation, that
is $\mathbb{G}_{n}(H)=O_{\mathbb{P}}(1).$ Formula (\ref{expan}) becomes%
\begin{equation*}
H_{n}=\mu (H)+n^{-1/2}\mathbb{G}_{n}(H)=\mu (H)+O_{\mathbb{P}}(n^{-1/2})
\end{equation*}

\bigskip \noindent and we will be able to use the Delta-method. Indeed, let $%
g:\mathbb{R}\longmapsto \mathbb{R}$ be continuously differentiable on a neighborhood of $\mu (H).$
The mean value theorem leads to%

\begin{equation}
g(H_{n})=g(\mu (H))+g^{\prime }(\mu _{n}(H))\text{ }n^{-1/2}\mathbb{G}_{n}(H),
\label{expan01}
\end{equation}

\bigskip \noindent where 
\begin{equation*}
\mu _{n}(H)\in \lbrack (\mu (H)+n^{-1/2}\mathbb{G}_{n}(H))\wedge \mu
(H),(\mu (H)+n^{-1/2}\mathbb{G}_{n}(H))\vee \mu (H)],
\end{equation*}

\bigskip \noindent so that 
\begin{equation*}
\left\vert \mu _{n}(H)-\mu (H)\right\vert \leq n^{-1/2}\mathbb{G}_{n}(H)=O_{%
\mathbb{P}}(n^{-1/2}).
\end{equation*}

\bigskip \noindent Then $\mu _{n}(H)$ converges to $\mu(H)$ in
probability (denoted as: $\mu _{n}(H)$ $\rightarrow _{\mathbb{P}}\mu (H)).$ But
the convergence in probability to a constant is equivalent to the weak
convergence. Then $\mu _{n}(H)$ $\rightsquigarrow \mu (H).$ Using again the
continuous mapping theorem, $g^{\prime }(\mu _{n}(H))\rightsquigarrow
g^{\prime }(\mu (H))$ which in turn yields $g^{\prime }(\mu
_{n}(H))\rightarrow _{\mathbb{P}}g^{\prime }(\mu (H))$ by the
characterization of the weak convergence to a constant. Now (\ref{expan01})
becomes

\begin{eqnarray*}
g(H_{n}) &=&g(\mu (H))+(g^{\prime }(\mu (H)+o_{P}(1))\text{ }n^{-1/2}\mathbb{%
G}_{n}(H) \\
&=&g(\mu(H))+ n^{-1/2} g^{\prime }(\mu (H))\mathbb{G}_{n}(H)+o_{\mathbb{P}}(n^{-1/2}) \\
&=&g(\mu (H))+ n^{-1/2}\mathbb{G}_{n}(g^{\prime}(\mu(H))H)+o_{P}(n^{-1/2}).
\end{eqnarray*}

\bigskip \noindent We arrive at the final expansion

\begin{equation}
g(H_{n})=g(\mu (H))+ n^{-1/2}\mathbb{G}_{n}(g^{\prime}(\mu(H))H)+o_{P}(n^{-1/2}).  \label{expanFinal}
\end{equation}

\Bin We just get a \textit{GAR} in the \textit{fep}, for some constant $C$,

\begin{equation}
J_n= C+ n^{-1/2}\mathbb{G}_{n}(h)+o_{P}(n^{-1/2}).  \label{expanFinal2}
\end{equation}

\Bin for $h^2$ being $\mathbb{P}_X$ integrable. A \textit{GAR} like (\ref{expanFinal2}), may undergo the Delta-method to give (under the needed conditions)

\begin{equation}
g(J_n)= g(C)+ n^{-1/2}\mathbb{G}_{n}(g^\prime(C)h)+o_{P}(n^{-1/2}).  \label{expanFinal3}
\end{equation}

\Bin which is still of the form of (\ref{expanFinal2}).\\

\Ni The general methods works as following: from exact \textit{GAR}s or \textit{GAR}s like (\ref{expanFinal2}), we get others \textit{GAR}s of the same form. Fortunately, many statistics are algebraic combinations of \textit{GAR}s in the form (\ref{expanFinal2}). To get the final \textit{GAR}, we have the tools of handling these algebraic combinations in the next theorem.

\begin{lemma} \label{lemma.tool.2}
Let ($A_{n})$ and ($B_{n})$ be two sequences of statistics expressed in the following \textit{GAR}s: real-valued functions$\ of$ $z\in S.$ Suppose that $A_{n}=A+n^{-1/2}\mathbb{G}_{n}(L)+o_{\mathbb{P}}(n^{-1/2})$ and $B_{n}=B+n^{-1/2}\mathbb{G}%
_{n}(H)+o_{\mathbb{P}}(n^{-1/2})$, where $A$, $B$, $L$ and $H$ are defined accordingly as above. We have the following \textit{GAR}s:

\begin{equation*}
A_{n}+B_{n}=A+B+n^{-1/2}\mathbb{G}_{n}(L+H)+o_{\mathbb{P}}(n^{-1/2}),
\end{equation*}

\bigskip \noindent

\begin{equation*}
A_{n}B_{n}=AB+n^{-1/2}\mathbb{G}_{n}(BL+AH) + o_{\mathbb{P}}(n^{-1/2}).
\end{equation*}

\bigskip \noindent and if $B\neq 0$,

\begin{equation*}
\frac{A_{n}}{B_{n}}=\frac{A}{B}+n^{-1/2}\mathbb{G}_{n}\biggr(\frac{1}{B}L-\frac{A}{%
B^{2}}H\biggr)+ o_{\mathbb{P}}(n^{-1/2}).
\end{equation*}
\end{lemma}

\bigskip \noindent By putting together all the described steps in a smart
way, the methodology will lead us to a final result of the form%
\begin{equation*}
T_{n}=t+n^{-1/2}\mathbb{G}_{n}(h)+o_{P}(n^{-1/2})
\end{equation*}%
which entails the weak convergence%
\begin{equation*}
\sqrt{n}(T_{n}-t)=\mathbb{G}_{n}(h)+o_{P}(1)\rightsquigarrow N(0,\Gamma
(h,h)).
\end{equation*}

\subsection{A non-trivial example}

\noindent We are going to illustrate our \textit{fep} tool on
the plug-in estimator of the linear correlation coefficient of the pair of random variables $(X,Y)$, none of the coordinates being degenerate, defined as
follows

\begin{equation*}
\rho =\frac{\sigma _{xy}}{\sigma _{x}^{2}\sigma _{y}^{2}},
\end{equation*}%
where%
\begin{equation*}
\mu _{x}=\int x\text{ }dP_{X}(x),\text{ }\mu _{y}=\int x\text{ }dP_{X}(x),%
\text{ }\sigma _{xy}=\int (x-\mu _{x})(y-\mu _{y})dP_{(X,Y)}(x,y).
\end{equation*}%
\begin{equation*}
\sigma _{x}^{2}=\int (x-\mu _{x})^{2}d\mathbb{P}_{X}(x),\text{ }\sigma _{y}^{2}=\int
(x-\mu _{x})(y-\mu _{y})d\mathbb{P}_{X}(y).
\end{equation*}%

\Bin We also dismiss the case the case where $\left\vert \rho \right\vert =1$,
for which one of $X$ and $Y$ is an affine function of the other, for example 
$X=aY+b$.\\

\bigskip \noindent It is clear that centering the variables $X$
and $Y$ and normalizing them by their standard deviations $\sigma _{x}$ and $%
\sigma _{y}$ does not change the correlation coefficient $\rho$. So we consider the coordinates as centered and normalized, i.e.

\begin{equation*}
\mu _{x}=\text{ }\mu _{y}=0,\text{ }\sigma _{x}=\sigma _{y}=1.
\end{equation*}

\bigskip \noindent However, we will let these coefficients appear with their names and we only use their particular values at the conclusion stage.\\

\noindent Let us consider the plug-in estimator of $\rho $. To this
end, let $(X_{1},Y_{1}),$ $(X_{2},Y_{2}),...$ be a sequence independent
observations of $(X,Y).$ For each $n\geq 1,$ the plug-in estimator is the following

\begin{equation*}
\rho _{n}=\left\{ \frac{1}{n}\sum_{i=1}^{n}(X_{i}-\overline{X})(Y_{i}-%
\overline{Y})\right\} \left\{ \frac{1}{n^{2}}\sum_{i=1}^{n}(X_{i}-\overline{X%
})^{2}\times \sum_{i=1}^{n}(Y_{i}-\overline{Y})^{2}\right\} ^{-1/2}.
\end{equation*}

\newpage
\bigskip \noindent We are going to give the asymptotic theory of $\rho _{n}$
as an estimator of $\rho$. Introduce the notation

\begin{equation*}
\mu _{(p,x),(q,y)}=E((X-\mu _{x})^{p}(Y-\mu _{y})^{q}),\mu _{4,x}=E(X-\mu
_{x})^{4}\text{, }\mu _{4,x}=E(X-\mu _{x})^{4}).
\end{equation*}

\bigskip \noindent Here is our main Theorem.

\begin{theorem}
\label{theo1}

\bigskip \noindent Suppose that neither of $X$ and $Y$ is degenerated and
both have finite fourth moments and that $X^{3}Y$ and $XY^{3}$ have finite
expectations. Then, as $n\rightarrow \infty ,$%
\begin{equation*}
\sqrt{n}(\rho _{n}-\rho )\rightsquigarrow N(0,\sigma ^{2}),
\end{equation*}
\bigskip \noindent where 
\begin{eqnarray*}
\sigma ^{2} &=&\sigma _{x}^{-2}\sigma _{y}^{-2}(1+\rho ^{2}/2)\mu
_{(2,x),(2,y)}+\rho ^{2}(\sigma _{x}^{-4}\mu _{4,x}+\sigma _{y}^{-4}\mu
_{4,y})/4 \\
&&-\rho (\sigma _{x}^{-3}\sigma _{y}^{-1}\mu _{(3,x),(1,y)}+\sigma
_{x}^{-1}\sigma _{y}^{-3}\mu _{(1,x),(3,y)}).
\end{eqnarray*}
\end{theorem}

\bigskip \noindent This result enables to test independence between $X$ and $%
Y$, or to test non linear correlation in the following sense.

\bigskip

\begin{theorem}
\label{theo2} Suppose that the assumptions of Theorem \ref{theo1} hold. Then%
\newline
\bigskip \noindent \textbf{(1)} If $X$ and $Y$ are not linearly correlated,
that is $\rho =0$, we have 
\begin{equation*}
\sqrt{n}\rho _{n}\rightsquigarrow N(0,\sigma _{1}^{2}),
\end{equation*}

\bigskip \noindent where 
\begin{equation*}
\sigma _{1}^{2}=\sigma _{x}^{-2}\sigma _{y}^{-2}\mu _{(2,x),(2,y)}.
\end{equation*}

\bigskip \noindent \textbf{(2)} If $X$ and $Y$ are independent, then $\rho
=0,$ and\bigskip 
\begin{equation*}
\sqrt{n}\rho _{n}\rightsquigarrow N(0,1)
\end{equation*}
\end{theorem}

\Bin \textbf{Proofs}. We are going to use the function empirical process based
on the observations $(X_{i},Y_{i}),i=1,2,...$ that are independent copies of 
$(X,Y)$. Write%
\begin{equation*}
\rho _{n}=\frac{\frac{1}{n}\sum_{i=1}^{n}X_{i}Y_{i}-\overline{X}\text{ }%
\overline{Y}}{\left\{ \frac{1}{n}\sum_{i=1}^{n}X_{i}^{2}-\overline{X}%
^{2}\right\} ^{1/2}\left\{ \frac{1}{n}\sum_{i=1}^{n}Y_{i}^{2}-\overline{Y}%
^{2}\right\} ^{1/2}}=\frac{A_{n}}{B_{n}}.
\end{equation*}

\bigskip \noindent Let us say for once that all the functions of $(X,Y)$
that will appear below are measurable and have finite second moments. Let us
handle separately the numerator and denominator. To treat $A_{n}$, using the
empirical process implies that

\begin{equation}
\left\{ 
\begin{tabular}{l}
$\frac{1}{n}\sum_{i=1}^{n}X_{i}Y_{i}=\mu _{xy}+n^{-1/2}G_{n}(p),$ \\ 
$\overline{X}=\mu _{x}+n^{-1/2}G_{n}(\pi _{1}),$ \\ 
$\overline{Y}=\mu _{y}+n^{-1/2}G_{n}(\pi _{2}),$%
\end{tabular}%
\right.  \label{casAN}
\end{equation}

\bigskip \noindent where $p(x,y)=xy$, $\pi _{1}(x,y)=x,$ $\pi _{2}(x,y)=y$.
From there we use the fact that $G_{n}(g)=O_{P}(1)$ for $\mathbb{E}%
(g(X,Y)^{2})<+\infty $ and get%
\begin{equation}
A_{n}=\mu _{xy}+n^{-1/2}G_{n}(p)-(\mu _{x}+n^{-1/2}G_{n}(\pi _{1}))(\mu
_{y}+n^{-1/2}G_{n}(\pi _{2})).  \label{haut}
\end{equation}
\bigskip \noindent This leads to 
\begin{equation*}
A_{n}=\sigma _{xy}+n^{-1/2}G_{n}(H_{1})+o_{P}(n^{-1/2})
\end{equation*}

\bigskip \noindent with 
\begin{equation*}
H_{1}(x,y)=p(x,y)-\mu _{x}\pi _{2}-\mu _{y}\pi _{1.}
\end{equation*}

\bigskip \noindent Next, we have to handle $B_{n}.$ Since the roles of $%
\left\{ \frac{1}{n}\sum_{i=1}^{n}X_{i}^{2}-\overline{X}^{2}\right\} ^{1/2}$
and $\left\{ \frac{1}{n}\sum_{i=1}^{n}Y_{i}^{2}-\overline{Y}^{2}\right\}
^{1/2}$ ar symmetrical, we treat one of them and extend the results to the
other. Consider $\left\{ \frac{1}{n}\sum_{i=1}^{n}X_{i}^{2}-\overline{X}%
^{2}\right\} ^{1/2}.$ From (\ref{casAN}), use the Delta-method to get%
\begin{equation*}
\overline{X}^{2}=\left( \mu _{x}+n^{-1/2}G_{n}(\pi _{1})\right) ^{2}=\mu
_{x}^{2}+2\mu _{x}n^{-1/2}G_{n}(\pi _{1})+o_{P}(n^{-1/2}),
\end{equation*}

\bigskip \noindent that is%
\begin{equation*}
\overline{X}^{2}=\left( \mu _{x}+n^{-1/2}G_{n}(\pi _{1})\right) ^{2}=\mu
_{x}^{2}+n^{-1/2}G_{n}(2\mu _{x}\pi _{1})+o_{P}(n^{-1/2}).
\end{equation*}

\bigskip \noindent From there, we get%
\begin{eqnarray*}
\frac{1}{n}\sum_{i=1}^{n}X_{i}^{2}-\overline{X}^{2}
&=&m_{2,x}+n^{-1/2}G_{n}(\pi _{1}^{2})-\overline{X}^{2} \\
&=&m_{2,x}-\mu _{x}^{2}+n^{-1/2}G_{n}(\pi _{1}^{2}-2\mu _{x}\pi
_{1})+o_{P}(n^{-1/2}) \\
&=&\sigma _{x}^{2}+n^{-1/2}G_{n}(\pi _{1}^{2}-2\mu _{x}\pi
_{1})+o_{P}(n^{-1/2}).
\end{eqnarray*}

\bigskip \noindent Using the Delta-method once again leads to 
\begin{equation*}
\left\{ \frac{1}{n}\sum_{i=1}^{n}X_{i}^{2}-\overline{X}^{2}\right\}
^{1/2}=\sigma _{x}+n^{-1/2}G_{n}\left(\frac{1}{2\sigma _{x}}\left\{ \pi
_{1}^{2}-2\mu _{x}\pi _{1}\right\} \right)+o_{P}(n^{-1/2}).
\end{equation*}

\bigskip \noindent In a similar way, we get%
\begin{equation*}
\left\{ \frac{1}{n}\sum_{i=1}^{n}Y_{i}^{2}-\overline{Y}^{2}\right\}
^{1/2}=\sigma _{y}+n^{-1/2}G_{n}\left(\frac{1}{2\sigma _{y}}\left\{ \pi
_{2}^{2}-2\mu _{y}\pi _{2}\right\} \right)+o_{P}(n^{-1/2}).
\end{equation*}

\bigskip \noindent We arrive at%
\begin{eqnarray*}
B_{n} &=&\left\{ \frac{1}{n}\sum_{i=1}^{n}X_{i}^{2}-\overline{X}^{2}\right\}
^{1/2}\left\{ \frac{1}{n}\sum_{i=1}^{n}Y_{i}^{2}-\overline{Y}^{2}\right\}
^{1/2} \\
&=&\sigma _{x}\sigma _{y}+n^{-1/2}G_{n}(\frac{\sigma _{y}}{2\sigma _{x}}%
\left\{ \pi _{1}^{2}-2\mu _{x}\pi _{1}\right\} +\frac{\sigma _{x}}{2\sigma
_{y}}\left\{ \pi _{2}^{2}-2\mu _{y}\pi _{2}\right\} )+o_{P}(n^{-1/2}).
\end{eqnarray*}

\bigskip \noindent By setting%
\begin{equation*}
H_{2}(x,y)=\frac{\sigma _{y}}{2\sigma _{x}}\left\{ \pi _{1}^{2}-2\mu _{x}\pi
_{1}\right\} +\frac{\sigma _{x}}{2\sigma _{y}}\left\{ \pi _{2}^{2}-2\mu
_{y}\pi _{2}\right\},
\end{equation*}

\bigskip \noindent we have%
\begin{equation}
B_{n}=\sigma _{x}\sigma _{y}+n^{-1/2}G_{n}(H_{2})+n^{-1/2}.  \label{bas}
\end{equation}

\bigskip \noindent Now, combining (\ref{haut}) and (\ref{bas}) and using
Lemma \ref{lemma.tool.1} yield%
\begin{equation*}
\sqrt{n}(\rho _{n}^{2}-\rho ^{2})=n^{-1/2}G_{n}\left(\frac{1}{\sigma _{x}\sigma
_{y}}H_{1}-\frac{\sigma _{xy}}{\sigma _{x}^{2}\sigma _{y}^{2}}%
H_{2}\right)+o_{P}(1).
\end{equation*}

\bigskip \noindent Let%
\begin{equation*}
H=\frac{1}{\sigma _{x}\sigma _{y}}(p(x,y)-\mu _{x}\pi _{2}-\mu _{y}\pi _{1})-%
\frac{\rho }{\sigma _{x}\sigma _{y}}\left\{ \frac{1}{2\sigma _{x}^{2}}%
\left\{ \pi _{1}^{2}-2\mu _{x}\pi _{1}\right\} +\frac{1}{2\sigma _{y}^{2}}%
\left\{ \pi _{2}^{2}-2\mu _{y}\pi _{2}\right\} \right\} .
\end{equation*}

\bigskip \noindent Now we continue with the centered and normalized case to
get%
\begin{equation*}
H(x,y)=p(x,y)-\frac{\rho }{2}(\pi _{1}^{2}+\pi _{2}^{2})
\end{equation*}

\bigskip \noindent and 
\begin{equation*}
H(X,Y)=XY-\frac{\rho }{2}(X^{2}+Y^{2}).
\end{equation*}

\bigskip \noindent Denote%
\begin{equation*}
\mu _{(p,x),(q,y)}=E((X-\mu _{x})^{p}(Y-\mu _{y})^{q}).
\end{equation*}

\bigskip \noindent We have 
\begin{equation*}
EH(X,Y)=\sigma _{xy}-\rho =0
\end{equation*}

\bigskip \noindent and $varH(X,Y)$ is equal to 
\begin{equation*}
\mu _{(2,x),(2,y)}+\rho ^{2}(\mu _{4,x}+\mu _{4,y})/4-\rho (\mu
_{(3,x),(1,y)}+\mu _{(1,x),(3,y)})+\rho ^{2}\mu _{(2,x),(2,y)}/2
\end{equation*}%
and finally $varH(X,Y)=\sigma _{0}^{2}$ with%
\begin{equation*}
\sigma _{0}^{2}=(1+\rho ^{2}/2)\mu _{(2,x),(2,y)}+\rho ^{2}(\mu _{4,x}+\mu
_{4,y})/4-\rho (\mu _{(3,x),(1,y)}+\mu _{(1,x),(3,y)}).
\end{equation*}

\bigskip \noindent This gives the conclusion that for centered and
normalized $X$ and $Y,$%
\begin{equation*}
\sqrt{n}(\rho _{n}-\rho )\rightsquigarrow N(0,\sigma _{0}^{2}).
\end{equation*}

\bigskip \noindent Next, if we use the normalizing coefficients in $\sigma
_{0}$, we get%
\begin{eqnarray*}
\sigma ^{2} &=&\sigma _{x}^{2}\sigma _{y}^{2}(1+\rho ^{2}/2)\mu
_{(2,x),(2,y)}+\rho ^{2}(\sigma _{x}^{4}\mu _{4,x}+\sigma _{y}^{4}\mu
_{4,y})/4 \\
&&-\rho (\sigma _{x}^{3}\sigma _{y}\mu _{(3,x),(1,y)}+\sigma _{x}\sigma
_{y}^{3}\mu _{(1,x),(3,y)})
\end{eqnarray*}

\bigskip \noindent and we conclude in the general case that%
\begin{equation*}
\sqrt{n}(\rho _{n}-\rho )\rightsquigarrow N(0,\sigma ^{2}).
\end{equation*}

\bigskip \noindent The proof of Theorem \ref{theo2} follows by easy computations under the particular conditions of $\rho $ and under
independence. $\blacksquare$\\

\newpage
\section{The Lo's residual functional emprical process \textit{(lrfep)} and its use in asymptotic theory} \label{sec_03}

\Ni Before we begin introducing the \textit{lrfep}, we need additional notions and notation around the R\'enyi's representation. Next we will explain how the \textit{lrfep} originated from the study ot the weak convergence of statistics. Finally, we will use it to get the full \textit{GAR} expression of statistics in general. 

\subsection{Additional notation} \label{sec_02_subsec_01} 

\Ni  In this part, we complete the notations we already gave and
specify our probability space.\newline

\noindent The notation below concerns the univariate random distributions. We are going to describe the
general Gaussian field in which we present our results. Indeed, we use a
unified approach when dealing with the asymptotic theories of the welfare
statistics. It is based on the Functional Empirical Process (\textit{fep})
and its Functional Brownian Bridge (\textit{fbb}) limit. It is laid out as
follows.\newline

\noindent When we deal with the asymptotic properties of one statistic or
index at a fixed time, we suppose that we have a non-negative random
variable of interest which may be the income or the expense $X$ whose
probability law on $(\mathbb{R},\mathcal{B}(\mathbb{R}))$, the Borel
measurable space on $\mathbb{R}$, is denoted by $\mathbb{P}_{X}.$ We
consider the space $\mathcal{F}_{(1)}$ of measurable real-valued functions $%
f $ defined on $\mathbb{R}$\ such that 

\begin{equation*}
V_{X}(f)=\int (f-\mathbb{E}_{X}(f))^{2}d\mathbb{P}_{X}=\mathbb{E}(f(X)-%
\mathbb{E}(f(X))^{2}<+\infty ,
\end{equation*}

\noindent where 
\begin{equation*}
\mathbb{E}_{X}(f)=\mathbb{E}f(X).
\end{equation*}

\bigskip \noindent On this functional space $\mathcal{F}_{(1)},$\ which is
endowed with the $L_{2} $-norm 
\begin{equation*}
\left\Vert f\right\Vert _{2}=\left( \int f^{2}d\mathbb{P}_{X}\right) ^{1/2},
\end{equation*}

\noindent we define the Gaussian process $\{\mathbb{G}_{(1)}(f),f\in 
\mathcal{F}_{(1)}\},$ which is characterized by its variance-covariance
function

\begin{equation}
\Gamma_{(1)}(f,g)=\int^{2}(f-\mathbb{E}_{X}(f))(g-\mathbb{E}_{X}(g))d\mathbb{P}_{X},(f,g)\in \mathcal{F}_{(1)}^{2}. \label{Gamma1}
\end{equation}

\noindent This Gaussian process is the asymptotic weak limit of the sequence of functional empirical processes (fep) defined as follows. Let $X_{1},X_{2},...$ be a sequence of independent copies of $X$. For each $n\geq1$, we define the functional empirical process associated with $X$ by 
\begin{equation*}
\mathbb{G}_{n,(1)}(f)=\frac{1}{\sqrt{n}}\sum_{j=1}^{n}(f(X_{j})-\mathbb{E}f(X_{j})),f\in \mathcal{F}_{(1)},
\end{equation*}

\noindent and denote the integration with respect to the empirical measure by

\begin{equation*}
\mathbb{P}_{n,(1)}(f)=\frac{1}{n} \sum_{i=1}^{n}f(X_{i}), \ f\in \mathcal{F}%
_{(1)},
\end{equation*}

\noindent Let us denote by $\ell^{\infty }(T)$ the space of real-valued bounded
functions defined on $T=\mathbb{R}$ equipped with its uniform topology. In
the terminology of the weak convergence theory, the sequence of objects $%
\mathbb{G}_{n,(1)}$ weakly converges to $\mathbb{G}_{(1)}$ in $\ell^{\infty}(%
\mathbb{R})$, as stochastic processes indexed by $\mathcal{F}_{(1)}$,
whenever it is a Donsker class. The details of this highly elaborated theory
may be found in \cite{billingsley68}, \cite{pollard},  \cite{vaart} and similar sources.\newline

\noindent Here, we only need the convergence in finite distributions which is a simple consequence of the multivariate central limit theorem, as described in Chapter 3 in \cite{ips-wcia-ang}.\newline

\noindent We will use the R\'enyi's representation of the random variable $X_i$'s of interest by means of the (\textit{cdf}) $F_{(1)}$ of  $X$ as follows 
\begin{equation*}
X=_{d}F_{(1)}^{-1}(U),
\end{equation*}

\noindent where $U$ is a uniform random variable on $(0,1)$, $=_{d}$ stands
for the equality in distribution and $F^{-1}_{(1)}$ is the generalized
inverse of $F_{(1)}$, defined by

\begin{equation*}
F_{(1)}^{-1}(s)=\inf \{x, F_{(1)}(x)\geq s\}, \ s \in (0,1).
\end{equation*}

\bigskip \noindent Based on these representations, we may and do assume that we are
on a probability space $(\Omega,\mathcal{A},\mathbb{P})$ holding a sequence
of independent $(0,1)$-uniform random variables $U_1$, $U_2$, ..., and the
sequence of independent observations of $X$ are given by

\begin{equation}
X_{1}=F_{(1)}^{-1}(U_1), \ \ X_{2}=F_{(1)}^{-1}(U_2), \ \ etc.
\label{repRenyi}
\end{equation}

\bigskip  \noindent For each $n\geq 1$, the order statistics of $U_1,...,U_n$
and of $X_1,...,X_n$ are denoted respectively by $0 \equiv U_{0,n} < U_{1,n}\leq \cdots \leq
U_{n,n} <U_{n+1,n}\equiv 1$ and $-\inf\equiv X_{0,n}<X_{1,n}\leq \cdots \leq X_{n,n} <X_{n+1,n}$ ($0=-\inf\equiv X_{0,n}$ for 
a non-negative random variable $X$).\newline

\noindent \noindent We also associate the sequence of $(U_n)_{n\geq 1}$ with the sequence of real-argument empirical functions

\begin{eqnarray}
\mathbb{U}_{n,(1)}(s)&=&\frac{1}{n} \#\{j,1\leq j \leq n, \ U_j \leq s\}, \
s\in(0,1) \ n\geq 1  \label{empiricalFunctionU}\\
&=& \sum_{j=1}^{n} \frac{j}{n} 1_{({U_{j,n}\leq s < U_{j+1,n}})},
\end{eqnarray}

\bigskip \noindent and the sequence of real-argument uniform quantile function

\begin{equation}
\mathbb{V}_{n,(1)}(s)=U_{1,n}1_{(s=0)}+\sum_{j=1}^{n}U_{j,n}1_{((j-1)/n < s \leq (j/n))}, \ s\in(0,1), \
n\geq 1  \label{quantileFunctionU}
\end{equation}

\bigskip \noindent and next, the sequence of real-argument uniform empirical processes 

\begin{equation}
\alpha_{n,(1)}(s)=\sqrt{n}(\mathbb{U}_{n,(1)}-s), 
\label{empiricalProcessU}
\end{equation}

\bigskip \noindent for $s\in (0,1)$ and $\geq 1$, and the sequence of real-argument uniform quantile processes

\begin{equation}
\gamma_{n,(1)}(s)=\sqrt{n}(s-\mathbb{V}_{n,(1)}), \ s\in(0,1) \ n\geq 1.
\label{quantileProcessU}
\end{equation}

\bigskip \noindent The same can be done for the sequence $(X_n)_{n\geq 1}$,
and we obtain the associated sequence of real empirical processes as

\begin{equation}
\mathbb{G}_{n,r,(1)}(x)=\sqrt{n} \left( \mathbb{F}_{n,(1)}(x)-F_{(1)}(x)%
\right), \ x\in \mathbb{R}, \ n\geq 1,  \label{empiricalProcess}
\end{equation}

\bigskip \noindent where

\begin{equation}
\mathbb{F}_{n,(1)}(x)=\frac{1}{n} \#\{j,1\leq j \leq n, \ X_j \leq x\}, \ x\in \mathbb{R} \ n\geq 1  \label{empiricalFunction}
\end{equation}

\noindent is the associated sequence of empirical functions. We also have the associated sequence of quantile processes

\begin{equation}
\mathbb{Q}_{n,(1)}(x)=\sqrt{n} \left( \mathbb{F}^{-1}_{(n),(1)}(s) -
F^{-1}(s) \right), \ s\in (0,1), \ n\geq 1  \label{quantileProcess}
\end{equation}

\noindent where, for $n\geq 1$,

\begin{equation}
\mathbb{F}^{-1}_{n,(1)}(s)=X_{1,n}1_{(0\leq s \leq 1/n)}+\sum_{j=1}^{n}X_{j,n}1_{((j-1)/n\leq s \leq (j/n))}, \ s\in(0,1),
\label{quantileFunction}
\end{equation}

\noindent is the associated sequence of quantile processes.\newline

\noindent By passing, we recall that $\mathbb{F}^{-1}_{n,(1)}$ is actually
the generalized inverse of $\mathbb{F}_{(n),(1)}$ and for the uniform
sequence, we have

\begin{equation}
\mathbb{V}_{n,(1)}=\mathbb{U}^{-1}_{n,(1)}.  \label{invCDFEMP}
\end{equation}

\noindent In virtue of Representation (\ref{repRenyi}), we have the
following remarkable relations

\begin{equation}
\mathbb{G}_{n,r,(1)}(x)=\alpha_{n,(1)}(F_{(1)}(x)), \ x\in \mathbb{R} \label{EmpEmpprocess}
\end{equation}

\bigskip \noindent and

\begin{equation}
\mathbb{Q}_{n,(1)}(x)=\sqrt{n}\left( F^{-1}_{(1)}(\mathbb{V}_{n,(1)}(s))-
F^{-1}_{(1)}(s)\right) \ s\in (0,1), \ n\geq 1.  \label{QQprocess}
\end{equation}

\bigskip \noindent We also have the following relations between the empirical functions and quantile functions

\begin{equation}
\mathbb{F}_{n,(1)}(x)=\mathbb{U}_{n,(1)}(F_{(1)}(x)), \ x\in \mathbb{R}
\label{EEFunction}
\end{equation}

\bigskip \noindent and

\begin{equation}
\mathbb{F}^{-1}_{n,(1)}(s)=F^{-1}_{(1)}(\mathbb{V}_{(n),(1)}(s)), \ s\in
(0,1), \ n\geq 1.  \label{QQFunction}
\end{equation}

\noindent As well, the real-argument and function-argument empirical processes are related as follows: for $n\geq 1$,

\begin{equation}
\mathbb{G}_{n,r,(1)}(x)=\mathbb{G}_{n,(1)}(f_{x}^{\ast}), \
\alpha_{n,(1)}(s)=\mathbb{G}_{n,(1)}(\tilde{f}_s), \ s \in (0,1), \ x \in \mathbb{R},
\label{empiricalProcessRealFonct}
\end{equation}

\bigskip \noindent where for any $x \in \mathbb{R}$, $f_{x}^{\ast}=1_{]-%
	\infty,x]}$ is the indicator function of $]-\infty,x]$ and for $s \in (0,1)$%
, $f_s=1_{[0,s]}$ and $\tilde{f}_s=1_{]-\infty,F^{-1}_{(1)}(s)]}$.\newline

\bigskip \noindent To finish the description, a result of Kiefer-Bahadur (See \cite{bahadur66}) that says that the addition of the sequences of uniform empirical processes and quantiles processes (\ref{empiricalProcessU}) and (\ref{quantileProcessU}) is asymptotically, and uniformly on $[0,1]$, zero in probability, that is

\begin{equation}
\sup_{s\in [0,1]} \left\vert \alpha_{n,(1)}(s)+\gamma_{n,(1)}(s) \right\vert
=o_{\mathbb{P}}(1) \text{ as } n\rightarrow +\infty.  \label{bahadurRep}
\end{equation}

\noindent This result is a powerful tool to handle the rank statistics when our studied statistics are $L$-statistics.\newline

\subsection{The origin of the lrfep} \label{sec_02_subsec_02} 

\noindent The method is developed here for a continuous \textit{cdf} $F_{(1)}$. There is a considerable class of statistics which are combinations of one dimensional statistics of the form

$$
L_n=d_n \sum_{1\leq j \leq n} c(j,n) q_0(X_{j,n}), \ n\geq 1,
$$

\bigskip \noindent where $q_0$ is some measurable mapping, $c(\circ,n)$ a function of $j \in \{1,\cdots,n\}$ and $(d_n)_{n\geq 1}$ is a sequence of real numbers. If $F_{(1)}$ is continuous, we may use the rank statistics $(R_{1,n}, \cdots, R_{n,n})$ defined by

$$
\forall 1\leq i \leq n, \ \forall 1\leq j \leq n, \ R_{j,n}=i \Leftrightarrow X_{i,n}=X_{j}.
$$

\bigskip \noindent But it happens that for any $n\geq 1$, for any $1\leq j \leq n$,

$$
\frac{R_{j,n}}{n}=\mathbb{F}_{n,(1)}(X_j),
$$

\bigskip \noindent and this leads to

$$
L_n=\frac{1}{n} \sum_{1\leq j \leq n} \biggr(n d_n c\left(n\mathbb{F}_{n,(1)}(X_j)\right)\biggr) q_0(X_{j}), \ n\geq 1.
$$ 

\bigskip \noindent The question is: how to treat such statistics? Obviously, this is a L-statistic and the theory of Chernoff might help. However, we will use another alternative developed over the years since that uses the so-called \textit{lrfep}. The method is based on the following facts: In many cases, there exists a measurable mapping $q_1$ such that 
$\mathbb{E}q_1(X)<+\infty$ and we may write

\begin{eqnarray*}
	L_n&=&\frac{1}{n} \sum_{1\leq j \leq n} q_1(F_{(1)}(X_j)) q_0(X_{j})\\
	&+&\frac{1}{n} \sum_{1\leq j \leq n} \biggr(n d_n c\left(n\mathbb{F}_{n,(1)}(X_j)\right)- q_1(F_{(1)}(X_j))\biggr) q_0(X_{j})\\
\end{eqnarray*}

\bigskip \noindent in a way that, by using the mean value theorem, we arrive at

\begin{eqnarray*}
	&&\frac{1}{n} \sum_{1\leq j \leq n} \biggr(n d_n c\left(n\mathbb{F}_{n,(1)}(X_j)\right)- q_1(F_{(1)}(X_j))\biggr) q_0(X_{j})\\
	&=&\frac{1}{n} \sum_{1\leq j \leq n} \biggr\{\mathbb{F}_{n,(1)}(X_j)-F_{(1)}(X_j)\biggr\}\ q_1(X_j) q_0(X_{j})+o_{\mathbb{P}}(n^{-1/2}),
\end{eqnarray*}

\bigskip \noindent where $q_1(\circ)$ is some measurable function, and further under specific conditions to be found and checked, we finally have the form, for $q=q_0q_1$,

$$
L_n=\frac{1}{n} \sum_{1\leq j \leq n} h(X_{j}) + \frac{1}{n} \sum_{1\leq j \leq n} \biggr\{\mathbb{F}_{n,(1)}(X_j)-F_{(1)}(X_j)\biggr\} \ q(X_j) + o_{\mathbb{P}}(n^{-1/2}).
$$  

\bigskip \noindent We conclude that, in our effort to asymptotically represent $L_n$ as an application of the empirical measure to some function $h$, that is $\mathbb{P}_n(h)$, we still have a residual term in the form of

$$
Re_n(\ell)=  \frac{1}{n}\sum_{j=1}^{n}\left(\mathbb{F}_{n,(1)}(X_j) - F_{(1)}(X_j) \right) q(X_j).
$$ 

\bigskip \noindent This made \cite{lrfep2010} to name it a residual empirical process and proceeded to its independent study. Now let us describe this stochastic process more deeply.

\subsection{The partial \textit{GAR} of \textit{lrfep} and conditions of its validity} \label{sec_02_subsec_03} 

\subsubsection{The main features of the \textit{lrfep}} \label{sec_02_subsec_03_01}

\Ni A residual empirical process is any stochastic process of the form

\begin{equation*}
Re_n(\ell)=  \frac{1}{n}\sum_{j=1}^{n}\left(\mathbb{F}_{n,(1)}(X_j) - F_{(1)}(X_j) \right) q(X_j).
\end{equation*} 

\bigskip \noindent where $q$ is a measurable function from $[0,1]$ to $\mathbb{R}$. Let us define

$$
\ell(s)=-q(F^{-1}_{(1)}(s)), \ s \in (0,1)
$$

\noindent and

$$
\Delta_n(q,s)=\biggr\{q\biggr(F_{(1)}^{-1}\left(\mathbb{V}_{n,(1)}(s)\right)\biggr)-q\biggr(F_{(1)}^{-1}(s)\biggr)\biggr\}, \ \ s\in (0,1). 
$$

\bigskip \noindent \textbf{We stress} that the function $\ell$ depends of the \textit{cdf} $F_{(1)}$ and should bave been denoted $\ell(\circ)=\ell(F_{(1)},\circ)$. This warning is important in the situation of spatial analysis, as we will see it.\\

\subsubsection{The general results on the \textbf{lrfep}} \label{sec_02_subsec_03_02}

\Ni We have the following results.

\begin{theorem} \label{theoGenRes} Let $F_{(1)}$ be continuous. If the following two assertions:\\
	
	\noindent (1) $\mathbb{E}q(X)<+\infty$ \label{cre1} \ \ (CRe1)\\
	
	\noindent and,\\
	
	\noindent (2) and, as $n \rightarrow +\infty$, \label{cre2}
	$$
	\int_{0}^{1} \sqrt{n} \left( s - \mathbb{V}_{n,(1)}(s)\right) \Delta_n(q,s)\, ds \rightarrow 0 \ \ (CRe2)
	$$ 
	
	\bigskip \noindent hold,  we have the representation
	
	$$
	\sqrt{n} Re_n(\ell) = \int_{0}^{1} \mathbb{G}_{n,(1)}(\tilde{f}_s)\,\ell(s)\,ds + o_p(1).
	$$
\end{theorem}

\noindent \textbf{Proof}. By using Formulas (\ref{empiricalFunctionU}) and (\ref{quantileFunctionU}), we get

\begin{equation*}
Re_n = \sum_{j=1}^{n} \int_{\frac{j-1}{n}}^{\frac{j}{n}}\left\{
\mathbb{F}_{n,(1)}(\mathbb{F}_{n,(1)}^{-1}(s)) - F_{(1)}(\mathbb{F}_{n,(1)}^{-1}(s))\right\}\, q\left(\mathbb{F}_{n,(1)}^{-1}(s)\right) \,ds,
\end{equation*}

\noindent and hence

\begin{equation}  \label{Rn}
Re_n = \int_{0}^{1}\left\{\mathbb{F}_{n,(1)}(F_{n,(1)}^{-1}(s)) -F_{(1)}(\mathbb{F}_{n,(1)}^{-1}(s))\right\}\, q\left(\mathbb{F}_{n,(1)}^{-1}(s)\right) \,ds.
\end{equation}

\bigskip

\noindent Before we proceed further, we will get some facts for Formulas (\ref{invCDFEMP}), (\ref{EmpEmpprocess}) and (\ref{QQprocess}). We already have

$$
\mathbb{F}_{n,(1)}(\circ)=\mathbb{U}_{n,(1)}\biggr(F_{(1)}(\circ)\biggr), \ \  and \ \ \mathbb{F}^{-1}_{n,(1)}(\circ)= F^{-1}_{(1)}\biggr(\mathbb{V}_{n,(1)}(\circ)\biggr),
$$

\Bin Next, we have (See page 148, \cite{ips-wcia-ang})

$$
F_{(1)}\biggr(F^{-1}_{(1)}(\circ))-0\biggr)\leq \circ \leq F_{(1)}\biggr(F^{-1}_{(1)}(\circ))+0\biggr).
$$

\Bin Next, by continuity of $F_{(1)}$, we have

$$
\mathbb{F}_{(1)}\biggr(\mathbb{F}^{-1}_{n,(1)}(\circ)\biggr)=\mathbb{V}_{n,(1)}(\circ)
$$

\Bin and 

$$
\mathbb{F}_{n,(1)}\biggr(\mathbb{F}^{-1}_{n,(1)}(\circ)\biggr)=\mathbb{U}_{n,(1)}\biggr(\mathbb{V}_{n,(1)}(\circ)\biggr).
$$

\Bin With this facts, we may handily transform $Re_n$ as follows:

\begin{eqnarray*}
	\sqrt{n} Re_n &=& \int_{0}^{1} \sqrt{n} \left\{\mathbb{U}_{n,(1)}\left( \mathbb{V}_{n,(1)}(s) \right) - \mathbb{V}_{n,(1)}(s)\right\}\, q\left(F_{(1)}^{-1}\left(\mathbb{V}_{n,(1)}(s)\right)\right) \, ds\\ 
	&=& \int_{0}^{1} \sqrt{n} \left( s - \mathbb{V}_{n,(1)}(s)\right) q\left(F_{(1)}^{-1}\left(\mathbb{V}_{n,(1)}(s)\right)\right) \, ds\\
	&+& \int_{0}^{1} \sqrt{n} \left( \mathbb{U}_{n,(1)}\left(\mathbb{V}_{n,(1)}(s)\right) -s\right)\, q\left(F_{(1)}^{-1}\left(\mathbb{V}_{n,(1)}(s)\right)\right)\, ds\\
	&= :& Re_n(1) + Re_n(2).
\end{eqnarray*}

\bigskip

\noindent From \cite{shwell} (page 585), we have

\begin{equation*}
\sup_{0\leq s\leq 1} \left| \mathbb{U}_{n,(1)}\left(\mathbb{V}_{n,(1)}(s) \right)
-s \right| \leq \frac{1}{n}.
\end{equation*}

\noindent We get

\begin{eqnarray}
	\left| Re_n(2) \right| &\leq&  \frac{1}{\sqrt{n}} \int_{0}^{1}  q\left(F_{(1)}^{-1}\left(\mathbb{V}_{n,(1)}(s)\right)\right)\, ds \notag\\
	&=& \frac{1}{\sqrt{n}} \left( \frac{1}{n} \sum_{j=1}^{n} q(X_j)\right), \ \ \label{CRe0}
\end{eqnarray}

\noindent which is an $o_{\mathbb{P}}(n^{-1/2})$ whenever $\mathbb{E}q(X)$ is finite. We obtain

\begin{eqnarray}
\sqrt{n} Re_n =  \int_{0}^{1} \sqrt{n} \left( s - \mathbb{V}_{n,(1)}(s)\right) q\biggr(F_{(1)}^{-1}\left(\mathbb{V}_{n,(1)}(s)\right)\biggr)\, ds + o_p(1).
\end{eqnarray}

\Ni Next, we try to replace $q\biggr(F_{(1)}^{-1}\left(\mathbb{V}_{n,(1)}(s)\right)\biggr)$ by $q\biggr(F_{(1)}^{-1}(s)\biggr)$. We will be led to

\begin{eqnarray}
&&\sqrt{n} Re_n\\
&&=  \int_{0}^{1} \sqrt{n} \left( s - \mathbb{V}_{n,(1)}(s)\right) q\biggr(F_{(1)}^{-1}(s)\biggr) \ ds \notag \\
&&+\int_{0}^{1} \sqrt{n} \left( s - \mathbb{V}_{n,(1)}(s)\right) \biggr\{q\biggr(F_{(1)}^{-1}\left(\mathbb{V}_{n,(1)}(s)\right)\biggr)-q\biggr(F_{(1)}^{-1}(s)\biggr)\biggr\}\, ds + o_p(1) \notag \\
&&=:\int_{0}^{1} \sqrt{n} \left( s - \mathbb{V}_{n,(1)}(s)\right) q\biggr(F_{(1)}^{-1}(s)\biggr)  \ ds \notag\\
&&+ \int_{0}^{1} \sqrt{n} \left( s - \mathbb{V}_{n,(1)}(s)\right) \Delta_n(q,s) \ ds+o_p(1).\notag 
\end{eqnarray}

\newpage
\Bin Finally by the the Bahadur's representation \eqref{bahadurRep} and by Condition (Re2), we finish the proof in the form:

\begin{eqnarray}
\sqrt{n} Re_n&=& \int_{0}^{1}\gamma_{n,(1)}(s) q\left(F^{-1}_{(1)}(s)\right) \,ds + o_p(1) \label{bahadur1}\\
&=&\int_{0}^{1}\gamma_{n,(1)}(s)+\alpha_{n,(1)} q\left(F^{-1}_{(1)}(s)\right) \,ds\notag \notag\\
&-& \int_{0}^{1} \alpha_{n,(1)}(s) q\left(F^{-1}_{(1)}(s)\right) \,ds + o_p(1) \notag\\
&=& -\int_{0}^{1} \biggr(\alpha_{n,(1)}(s)\biggr) q\left(F^{-1}_{(1)}(s)\right) \,ds + o_p(1) \ \ \ \notag\\
&=& \int_{0}^{1} \mathbb{G}_{n,(1)}(\tilde{f}_s)\, \left(-q\left(F^{-1}_{(1)}(s)\right)\right) \,ds + o_p(1), \notag
\end{eqnarray}

\Bin i.e,

\begin{equation*}
\sqrt{n} Re_n = \int_{0}^{1} \mathbb{G}_{n,(1)}(\tilde{f}_s)\,\ell(s)\,ds + o_p(1),
\end{equation*}

\noindent whenever 

$$
\mathbb{E}(\ell(X))=-\mathbb{E}q(X)=-\int_{0}^{1} q(F_{(1)}^{-1}(s)) \,ds<+\infty.
$$

\noindent This concludes the proof.\\

\subsection{The final \textit{GAR} combining the \textbf{fep} and the \textbf{lrfep}} \label{sec_02_subsec_03} 

 \Ni \textbf{Summary of the general method}.\\

\Ni We suppose that we study the sequence of statistics $J_n$. We split the study into steps.\\

\Ni \textbf{Step 1}. Proceed to stochastic calculus and manipulations to try to write as:

\begin{equation}
J_n= \mathbb{E}h(X) + n^{-1/2} \mathbb{G}_n(h)+ Re_n(\ell) + o_{\mathbb{P}}(n^{-1/2}), \label{mg_01}
\end{equation}

\Bin where $\ell(\circ)$ is derived from some $q(x)$ measurable function of $x \in \mathbb{R}$ by $\ell(s)=-q(F^{-1}_{(1)}(s))$, for $(0,1)$, and

\begin{equation}
Re_n(\ell)=\frac{1}{n} \sum_{1}^{n} \{\mathbb{F}_{(n),(1)}(X_j)-\mathbb{F}_{(1)}(X_j)\} \ q(X_j), \ \ n\geq 1. \label{mg_02}
\end{equation}

\Bin We recall that

\begin{equation}
\Delta_n(q,s)=\biggr\{q\biggr(F_{(1)}^{-1}\left(\mathbb{V}_{n,(1)}(s)\right)\biggr)-q\biggr(F_{(1)}^{-1}(s)\biggr)\biggr\}, \ \ s\in (0,1). 
\label{mg_03}
\end{equation}

\Ni \textbf{Step 2}. Using Theorem \ref{theoGenRes}, we try to get
	
\begin{equation}
\mathbb{E}q(X)<+\infty \label{mg_04}
\end{equation}
	
	\noindent and,\\
	
	\noindent (2) as $n \rightarrow +\infty$, \
	
	\begin{equation}
		\int_{0}^{1} \sqrt{n} \left( s - \mathbb{V}_{n,(1)}(s)\right) \Delta_n(s)\, ds \rightarrow 0 \label{mg_05}
	\end{equation}

\Bin If so, we get the \textit{GAR} that will be displayed below.\\
	
\Ni \textbf{Step 2 (alternative)} If we are able to prove that

\begin{equation}
\frac{1}{\sqrt{n}} \left( \frac{1}{n} \sum_{j=1}^{n} q(X_j)\right)=o_P(1), \ as \ n\rightarrow +\infty, \label{mg_06}
\end{equation}

\Bin and

\begin{equation}
\int_{0}^{1}  (s(1-s))^{1/2} \Delta_n(s)=o_P(1), \ as \ n\rightarrow +\infty, \ \ \label{mg_07}
\end{equation}

\Bin hold, then the general \textit{GAR} is true. That \textit{GAR} is:

\begin{equation}
J_n= \mathbb{E}h(X) + n^{-1/2} \mathbb{G}_n(h)+ n^{-1/2} \int_{0}^{1} \mathbb{G}_n(\tilde f_s) \ \ell(s) \ ds +o_{\mathbb{P}}(n^{-1/2}), \label{mg_08}
\end{equation}

\Bin with 

\begin{equation}
\tilde f_s=1_{]-\infty, F^{-1}_{(1)}(s)]}, \ s \ (0,1). \label{mg_09}
\end{equation}

\begin{equation}
\sqrt{n}(J_n- \mathbb{E}h(X))= \mathbb{G}_n(h)+  \int_{0}^{1} \mathbb{G}_n(\tilde f_s) \ \ell(s) \ ds +o_{\mathbb{P}}(1). \label{mg_10}
\end{equation}

\Ni To understand Condition \eqref{mg_07}, we used the 

$$
\sup_{s \in (0,1)} \|\frac{\sqrt{n} \left( s - \mathbb{V}_{n,(1)}(s)\right)}{(s(1-s))^{1/2}}\|=O_{\mathbb{P}}(1).
$$

\Ni Now, we may conclude. The finite-distribution asymptotic law of the field $\mathbb{G}_n(f), \ f\in L^2(\mathbb{P}_X)$ gives the limit law of $\sqrt{n}(J_n- \mathbb{E}h(X))$, as follows.

\begin{theorem}\label{theoGenFull} If $h\in L^{2}(\mathbb{P}_X)$ and, (\ref{mg_04}) or (\ref{mg_05}) and (\ref{mg_06}) or (\ref{mg_07}) hold, then the sequence of indexes $J_n$ as in (\ref{mg_08}) satisfies

$$
\sqrt{n}(J_n-\mathbb{E}h(X)) \rightsquigarrow \mathcal{N}(0, \ \Gamma^{(J)}), \ \ as \ \ n\rightarrow +\infty,
$$

\Bin where

\begin{equation}
\Gamma_{(1)}(h_1,h_2)=\int (h_1-\mathbb{E}h_1(X))(h_2-\mathbb{E}h_2(X)) \ d\mathbb{P}_X \label{varTheo1}
\end{equation}

\Bin and

\begin{eqnarray}
&&\Gamma^{(J)}=\gamma_{1} + \gamma_{2} + 2\gamma_{3}, \label{varTheo2}\\
&&\gamma_{1}= \Gamma_{(1)}(h,h); \notag\\
&&\gamma_{2}= \int^{1}_{0}\int^{1}_{0}\Gamma_{(1)}(\tilde f_s, \tilde f_t) \ell(s) \ \ell(t) \ ds \ dt; \notag\\
&&\gamma_{3}=\int^{1}_{0}\Gamma_{(1)}(h, \tilde f_s) \ell(s) \  ds. \notag
\end{eqnarray}
\end{theorem}

\Bin \textbf{(R) Important Remarks on the existence of the variances and co-variance in \eqref{varTheo1} and \eqref{varTheo2}}.\\

\Ni \textbf{(A1) Asymptotic law of} \label{(Th22)}

$$
Re_n(\ell)=\int^{1}_{0} \mathbb{G}_n(\tilde f_s) \ \ell(s) \ ds.
$$

\Bin Since $\tilde f_s$ and $\tilde f_t$ are in $[0,1]$, that asymptotic variance of  $Re_n(\ell)$, is finite 

$$
\mathbb{E}(q(X)^2)<+\infty. \ \ \ (Th33)
$$ 

\Ni \textbf{(A2) Asymptotic covariance between $\mathbb{G}_n(h)$ and $Re_n(\ell)$} \label{(Th12)}. For the same reason, that Asymptotic covariance is finite if

$$
\mathbb{E}(h(X)^2 q(X)^2)<+\infty. \ \ \ (Th12)
$$

\Ni \textbf{(A3) Asymptotic variance of $J_n$} \label{(Th11)}. Finally, that asymptotic variance is finite
if (Th22), (Th12) and

$$
\mathbb{E}h(X)^2 <+\infty. \ \ \ (Th11)
$$

\newpage
\section{Finding the \textbf{GAR} in Gini's Case study} \label{sec_04}

\Ni Before we begin, we want to draw the attention of the reader who wants to work on this topic about the following warning.\\

\Ni \textbf{Warning}. The method we are presenting here uses a lot of stochastic expansions over quantities as $o_{\mathbb{P}}(1)$ and 
$O_{\mathbb{P}}(1)$. The researcher should have to try to get a good command in such computation skills. He may need to go to the bottom of that by revising Chapter 6 in \cite{ips-wcia-ang}. Here, in the present study, we do not dwell in the details.\\ 

\Ni Now, let us consider the Gini index of inequality as in \cite{mergane2018Gini} (where the multiplicative coefficient (2) has been corrected to (1))

\begin{equation*} \label{nap}
J_n=\frac{1}{\frac{1}{n}\sum_{j=1}^{n}X_j} \biggr\{\frac{1}{n} \sum_{j=1}^{n} \biggr(\frac{2j-1}{n}-1\biggr)X_{j,n}\biggr\},
\end{equation*}

\Bin for $n\geq 1$ and the notations and terminology given above are adopted. To apply the method we expose above, we see that we have to treat the two factors

\begin{equation*} \label{nap}
A_n=\frac{1}{\frac{1}{n}\sum_{j=1}^{n}X_j} \ \ and \ \ B_n=\frac{1}{n} \sum_{j=1}^{n} \biggr(\frac{2j-1}{n}-1\biggr)X_{j,n}
\end{equation*}

\Bin separately. The first step is easy to handle, for $h_1(x)=x$ and by denoting $\mathbb{E}h_1(X)=\mathbb{E}X=\mu$,

\begin{eqnarray*} \label{nap}
A_n&=&\frac{1}{\frac{1}{n}\sum_{j=1}^{n}X_j}\\
&=&\frac{1}{\mathbb{E}h_1(X)+ n^{-1/2}\mathbb{G}_n(h_1)}\\
&=&\frac{1}{\mu} \biggr(1 + \mu^{-1} n^{-1/2}\mathbb{G}_n(h_1)\biggr)^{-1}\\
&=&\frac{1}{\mu} \biggr(1 - \mu^{-1} n^{-1/2}\mathbb{G}_n(h_1) + O_{\mathbb{P}}(n^{-1/2})\biggr).
\end{eqnarray*}

\Bin We get a first \textit{GAR}

\begin{equation*}
A_n=\frac{1}{\mu} -  n^{-1/2}\mathbb{G}_n(\mu^{-2} h_1) + o_{\mathbb{P}}(n^{-1/2}). 
\end{equation*}

\Bin As to $B_n$, we remark that

\begin{eqnarray*} \label{nap}
B_n&=&\biggr( \frac{1}{n}\sum_{j=1}^{n} \biggr(\frac{2j}{n}-1\biggr) X_{j,n}\biggr) - \biggr(\frac{1}{n^2}\sum_{j=1}^{n} X_j\biggr)\\
&=:&B_{n,1} + B_{n,2}. 
\end{eqnarray*}

\Bin We have

\begin{equation} \label{nap01}
\sqrt{n}B_{n,2} =n^{-1} (\mu + o_{\mathbb{P}}(1))=o_{\mathbb{P}}(n^{-1/2}).
\end{equation}

\Bin Next,

\begin{eqnarray*} \label{nap}
B_{n,1}&=&\frac{1}{n}\sum_{j=1}^{n} \biggr(\frac{2j}{n}-1\biggr) X_{j,n}\biggr)\\ 
&=&\frac{1}{n}\sum_{j=1}^{n} \biggr(2F_n(X_{j,n})-1\biggr) X_{j,n}\\
&=&\frac{1}{n}\sum_{j=1}^{n} \biggr(2F_n(X_j)-1\biggr) X_j\\
&=&\frac{1}{n}\sum_{j=1}^{n} \biggr(2F(X_j)-1\biggr) X_j\\
&+&\frac{1}{n}\sum_{j=1}^{n} \biggr(F_n(X_j)-F(X_j)\biggr) 2X_j\\
&=&\mathbb{P}_n(H)+ Re_n(\ell_1)
\end{eqnarray*}

\Bin with

\begin{equation*} \label{nap}
H(x)=(2F(x)-1)x \ \ \ and \ \ \ \ell_1(x)=-2F^{-1}(s), \ \ x \in\mathbb{R}, \ \ s \in (0,1).
\end{equation*}

\Bin So, upon checking Conditions (\ref{mg_04}) or (\ref{mg_05}) and (\ref{mg_06}) or (\ref{mg_07}), we would get

\begin{equation} \label{nap10}
B_{n}=\mathbb{E}H(X) + n^{-1/2} \mathbb{G}_n(H) + n^{-1/2} \int_{0}^{1} \mathbb{G}_n(\tilde f_s) \ \ell_1(s) \ ds+o_{\mathbb{P}}(n^{-1/2})
\end{equation}

\Bin Finally, algebraic manipulations of the product of the expansions of $A_n$ and $B_n$, based on $\mathbb{G}_n(\circ)=o_{\mathbb{P}}(1)$ anywhere here, lead to

\begin{equation*} \label{nap20}
J_n=\frac{\mathbb{E}H(X)}{\mu} + n^{-1/2} \mathbb{G}_n\left(\frac{H}{\mu}-\frac{\mathbb{E}H(X)h_1}{\mu^2}\right) + n^{-1/2} \int_{0}^{1}(1/\mu)\mathbb{G}_n(\tilde f_s) \  \ell_1(s) \ ds+o_{\mathbb{P}}(n^{-1/2}),
\end{equation*}

\Bin which gives us our final \textit{GAR}

\begin{equation} \label{nap20}
J_n=\frac{\mathbb{E}H(X)}{\mu} + n^{-1/2} \mathbb{G}_n(h) + n^{-1/2}  \int_{0}^{1} \mathbb{G}_n(\tilde f_s) \ell(s) \ ds+o_{\mathbb{P}}(n^{-1/2}),
\end{equation}

\Bin where

\begin{eqnarray} \label{nap}
&&h(x)=\frac{\mu H(x)-\mathbb{E}H(X)h_1(x)}{\mu^2},\\
&&h_1(x)=x,\\
&&H(x)=(2F(x)-1)x \label{functs}, \notag\\
&&\ell(s)=\frac{-2F^{-1}(s)}{\mu}, \ \ x \in\mathbb{R}, \ \ s \in (0,1). \notag
\end{eqnarray}

\section{Conclusion, perspectives and applications} \label{sec_05}

\noindent The theory reaches two kinds of results.\\

\Ni (R1) A definitive result of asymptotic normality as in Theorem \ref{theoGenFull} that can be applied in statistical estimation by confidence intervals and in statistical tests.\\

\Ni (R2) A progressive result in the form of the asymptotic representation as in (\ref{mg_08}) (page \pageref{mg_08}) can then be used for comparison with other statistics or with values of the same statistic for other areas, or for following the evolution of the values of
the statistics other times in a way to ensure convergence of joint laws.\\

\Ni The following remarks are important to be noticed.\\

\Ni (R3) Let us be aware that the final result may come after algebraic manipulations of a number of \textit{GAR} founded separately.\\

\Ni (R4) The computations of variances and co-variances of the kind of statistics used here are well-handled with current computers and numerical methods with weak execution times.\\

\Ni (R5) To be completely useful to researchers and practitioners, a package should be available and ready-to-be used. Such a work is underway and will be reported on soon. 



\begin{thebibliography}{99}
\bibitem[Bahadur (1966)]{bahadur66} Bahadur R.R. (1966). A note on quantiles in large samples, \textit{Ann. Math. Stat.} \textbf{37}, pp. 577--580. MR :32:6522. ZL: 0147.18805. doi:10.1214/aoms/1177699450. euclid.aoms/1177699450.
\bibitem[Billingsley (1968)]{billinsgley68} Billingsley P. (1968). \textit{Convergence of probability measures}. Wiley, new-York.
\bibitem[Cs\H{o}rg\"o \textit{et al.} (1986)]{cscshm} Cs\H{o}rg\H{o}, H., Cs\H{o}r\H{o}, M., Mason D.M. and Horv\`ath, L. (1986) Weighted empirical process and quantile process. \textit{Ann. Probab.} 14 (1), 31-85.
\bibitem[Foster \textit{et al.} (1984)]{fgt} Foster, J., Greer J. and Shorrocks A. (1984). A class of Decomposable Poverty Measures.
\textit{Econometrica} 52, 761-766.
\bibitem[Gaenssler (1983)]{gaenssler} Gaenssler, P. (1983). \textit{Empirical processes}. IMS Lecture Notes-Monograph Series, Hayward, California.
\bibitem[Haidara \textit{et al.} (2018)]{haidaraDecomp} Haidara M.C, Kpanzou T.A, Mergane P.D. and Lo G.S. (2018). Asymptotic Theory and Statistical Decomposability gap Estimation for Takayamas Index. In A Collection of Papers in Mathematics and Related Sciences, a festschrift in honour of the late Galaye Dia (Editors : Seydi H., Lo G.S. and Diakhaby A.). Spas Editions, Euclid Series Book, pp. 195 –220.
Doi : 10.16929/sbs/2018.100-02-08.
\bibitem[Jarque and Bera (1987)]{jarque} Jarque C.M. and Bera A.K. (1987). A test for normality of observations and regression residuals. \textit{International Statistical Review}. 55(5): 163-172.
\bibitem[Kakwani (1980)]{kakwani} Kakwani N. (1980). On a Class of Poverty Measures. \textit{Econometrica}, 48, 437-446.
\bibitem[Lo (2013)]{logpi} Lo G.S. (2013). On the General Poverty Index. \textit{Far East Journal of Theoretical Statistics}. Volume 42. (1), 1-22.
\bibitem[Lo \textit{et al.} (2016)]{ips-wcia-ang} Lo, G.S., Ngom M. and Kpanzou T.A. (2016). Weak Convergence (IA). Sequences of random vectors. SPAS Books Series. (2016). Doi : 10.16929/sbs/2016.0001. Arxiv : 1610.05415.
\bibitem[Lo and Sall (2010a)]{lrfep2010} Lo G.S. and Sall S.T. (2010a).  Asymptotic Representation Theorems for Poverty Indices. \textit{Afrika Statistika, Special Volume (5) : Proceedings of the International Workshop on Multiple Risks and Copula, Biskra 2010, pp. 238-244. Ed. Abdelhakim Necir}. 
\bibitem[Lo and Sall (2010b)]{lrfep2010TD} Lo G.S. and Sall S.T. (2010b). Uniform weak convergence of the time-dependent poverty measures for continuous longitudinal data, \textit{Brazilian Journal of Probability and Statistics}, Vol. 24, Issue 3, pp 457-467.
\bibitem[Lo et \textit{al.}(2015)]{gsloJB} Lo G.S., Thiam, O. and Haidara C.M. (2015) High Moments Jarque-Bera Tests for Arbitrary Distribution Functions. \textit{Applied Mathematics}. Vol. 6 (5), 707-716. http://dx.doi.org/10.4236/am.2015.64066.
\bibitem[Mergane \textit{et al.} (2018a)]{mergane2018Theil} Mergane  P.D., Kpanzou T.A., Ba D. and Lo G.S. (2018a). A Theil-like class of inequality measures, its asymptotic normality theory and applications. \textit{Afrika Statistika} Vol. 13 (3), pp 1699-1715.
\bibitem[Mergane \textit{et al.} (2018b)]{mergane2018Gini} Mergane P.D., Lo G.S. and Kpanzou T.A. (2018b). On the joint distribution of variations of the Gini index and Welfare indices. In A Collection of Papers in Mathematics and Related Sciences, a festschrift in honour of the late Galaye Dia. Spas Editions, Euclid Series Book, pp. 163– 194. Doi : 10.16929/sbs/2018.100-02-07.
\bibitem[Pollard (1984)]{pollard} Pollard D. (1984). \textit{Convergence of Stochastic Processes}. Springer-Verlag, Berlin.
\bibitem[Sen (1976)]{sen} Sen A.K. (1976). Poverty: An Ordinal Approach to Measurement. \textit{Econometrica}, 44, 219-231.
\bibitem[Shorack and Wellner (1986)]{shwell} Shorack G.R. and Wellner J.A. (1986). Empirical Processes with Applications to Statistics, wiley-Interscience, New-York.
\bibitem[van der Vaart and Wellner (1996)]{vaart} van der Vaart A.W. and Wellner J.A. (1996). \textit{Weak Convergence and Empirical Processes With Applications to Statistics}. Springer, New-York.
\bibitem[Vapnik and Chervonenkis (1971)]{vcclass} Vapnik V.N. and Chervonenkis A.Y. (1971). On the uniform convergence of relative frequencies of event to their probabilities. \textit{Theory of Probability and Its applications}. 16, 264-280.
\bibitem[Zeng (1997)]{zeng} Zheng B. (1997). Aggregate poverty measures. \textit{Journal of Economic Surveys}, 11 123–162.
\end{thebibliography}
\end{document}